\documentclass[12pt]{article}
\usepackage{graphicx}

\font\tendb=msbm10 at 12pt \font\sevendb=msbm10 at 9pt
\font\fivedb=msbm10 at 7pt
\newfam\dbfam
\textfont\dbfam=\tendb \scriptfont\dbfam=\sevendb
\scriptscriptfont\dbfam=\fivedb
\def\db{\fam\dbfam\tendb}
\font\eufm=eufm10\font\eufms=eufm10\font\eufmss=eufm10\newfam\eufam
\textfont\eufam=\eufm\scriptfont\eufam=\eufms\scriptscriptfont\eufam=\eufmss

\font\tendbb=msbm10 at 12pt \font\sevendbb=msbm7 at 9pt
\font\fivedbb=msbm5 at 6pt
\newfam\dbbfam
\textfont\dbbfam=\tendbb \scriptfont\dbbfam=\sevendbb
\scriptscriptfont\dbbfam=\fivedbb

 \def \Z {{\db Z}}
 \def \R {\hbox{\db R}}

\font\tenMmm=eusm10 at 12pt
\newfam\Mmmfam
\textfont\Mmmfam=\tenMmm

\def\illu #1 by #2 (#3){
  \vbox to #2{
    \hrule width #1 height 0pt depth 0pt
    \vfill
    \special{illustration #3} 
    }
  }

\begin{document}

\null

\vspace{2cm}

\begin{center}
{\large {\bf Ribbon graphs   and the Temperley-Lieb Algebra}}\\
 Nafaa Chbili\\
\begin{footnotesize} Department of Mathematical Sciences\\
 College of Science\\
 UAE University\\
 E-mail: nafaachbili@uaeu.ac.ae\\
 \end{footnotesize}

\end{center}


\begin{abstract}
Let $n$ be a nonnegative  integer, we use ribbon $n-$graph diagrams  and the Yamada polynomial skein relations
to construct an algebra ${\mathcal Y}_n$ which is shown  to be  closely related to the Temerley-Lieb Algebra. We prove that the algebra
 ${\mathcal Y}_2$ is isomorphic to some  quotient of a  three variables polynomial algebra. Then,
 we give a family of generators for the algebra  ${\mathcal Y}_3$.
\end{abstract}

\section{Introduction}
Throughout this paper, a graph is the geometric realization of a finite CW-complex of dimension 1. Furthermore, we assume that all vertices (0-cells)
have valency greater than 2.
A spatial graph is an embedding of a  graph into the three-dimensional Euclidean space $\R^3$. The theory of spatial graphs is considered  as  a natural  extension of knot theory. Therefore, many of the techniques and problems  of  knot theory have their counterparts in spatial graph theory. A natural question that arose  after the discovery of the Jones polynomial  and other quantum invariants of links, was to define   invariants of Jones type for  spatial graphs. In that direction, Yamada \cite{Ya1} introduced a topological  invariant of spatial  graphs, hereafter referred to as the Yamada polynomial. It is a one variable Laurent polynomial $Y(A)$ which can be defined recursively on planar diagrams of spatial graphs.\\
The   Jones and the Kauffman bracket polynomials are closely related to the Temperley-Lieb algebra ${\mathcal \tau}_n$. Actually, one can construct these polynomials through representations of the Artin braid group into ${\mathcal \tau}_n$. The original  motivation of the present paper is to
explore the possibility of a similar algebraic interpretation  of the Yamada polynomial.\\

Ribbon graphs are geometrical objects that  appeared as a natural generalization of framed links by  Reshetikhin and Turaev in \cite{RT}.  Let $n$ be a nonnegative  integer. A ribbon $n-$graph is a compact oriented surface embedded into $\R^2 \times I$ which meets the boundary of  $\R^2 \times I$  orthogonally exactly along the  $2n$ segments  $\{[i-1/10,i+1/10]\times\{0\}\times \{0,1\},i=1,\dots, n\}$. It is worth mentioning that the precise definition requires  some other technical arrangements the discussion of  which is postponed to Section 3. Ribbon graphs are  represented by planar graph diagrams generalizing link diagrams.  Let ${\mathcal S}_n$ be the set of all ribbon $n-$graph diagrams and  $\mathcal
R$=$\Z[A^{\pm 1},d^{-1}]$, where $d=-A^2-A^{-2}$. Let
${Y}_n$ be the free $\mathcal R$-module generated by all elements
of  ${\mathcal S}_n$. We define $\mathcal{Y}_n$ to  be the quotient module  of ${Y}_n$ by the Yamada relations in Section 3. This module admits a natural algebra structure. The product of two ribbon $n-$graphs $G$ and $G'$ is defined as   illustrated by  Figure 1.

\begin{center}
\includegraphics[width=5cm,height=4cm]{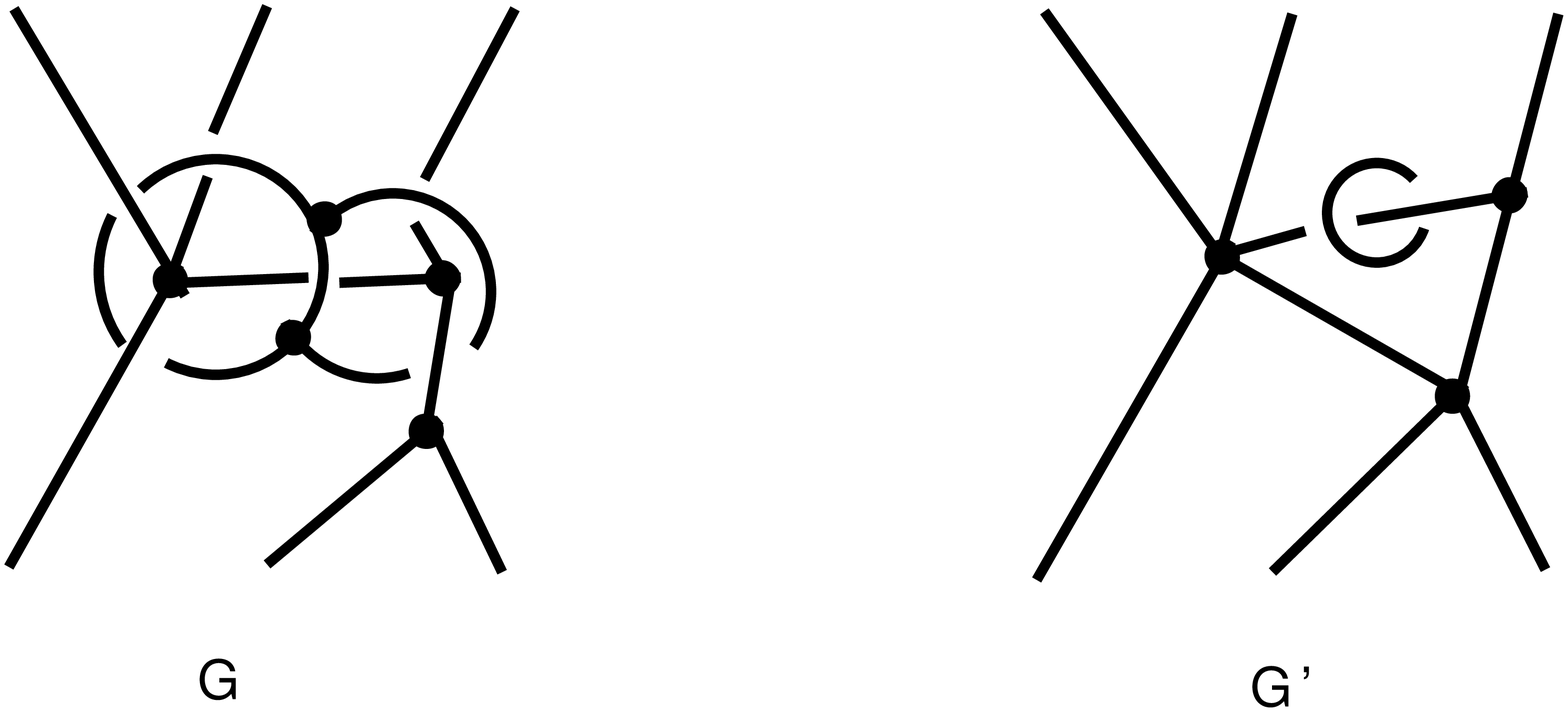} \hspace{2cm}
\includegraphics[width=2cm,height=4cm]{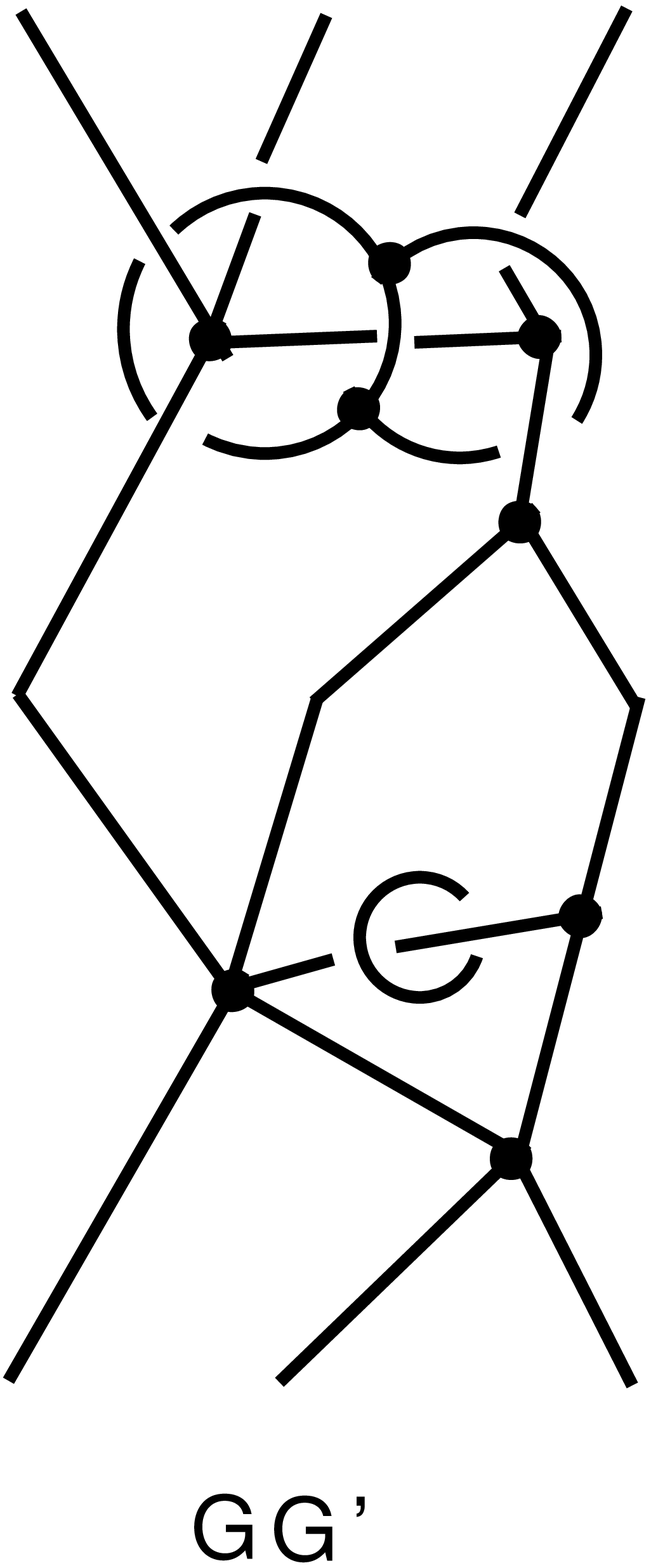}\\
Figure 1
\end{center}
Here are the main results in this paper.\\
\textbf{Theorem 1.1.} {\sl The algebra ${\mathcal Y}_2$ is the free additive $\mathcal R$-algebra with multiplicative elements $1_2, V$ and $X$ pictured below.}\
\begin{center}
\includegraphics[width=9cm,height=3cm]{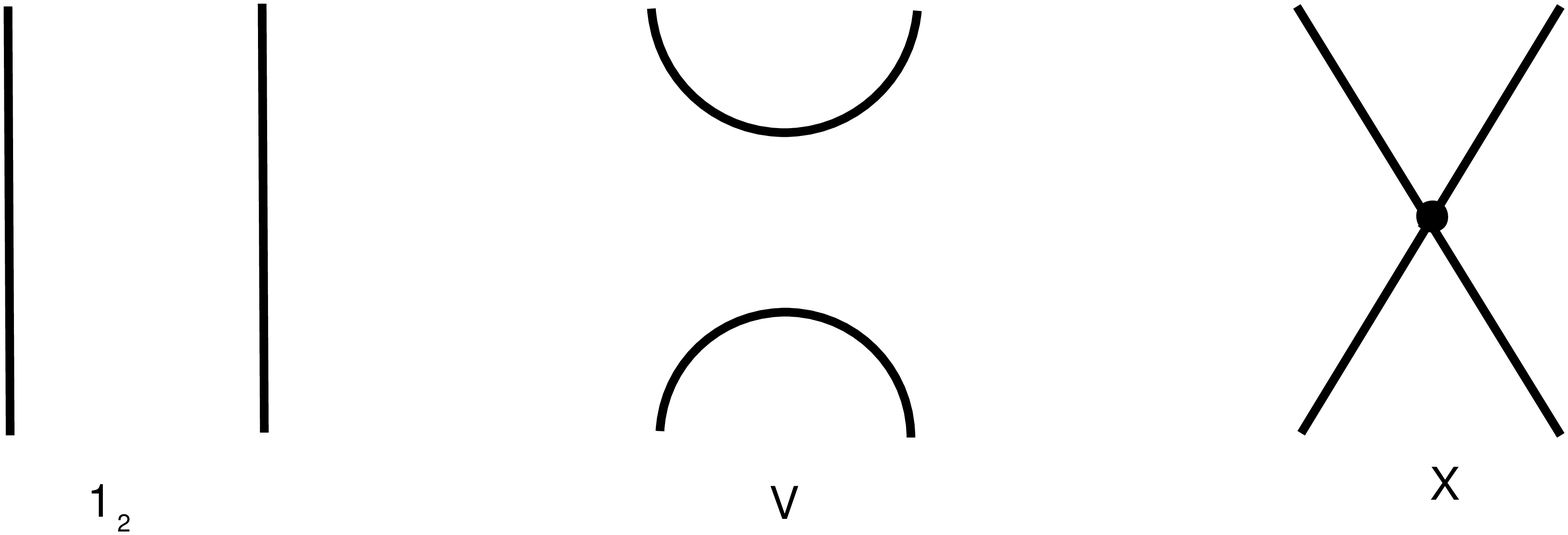} \\
Figure 2
\end{center}

\textbf{Corollary 1.1.} {\sl The algebra ${\mathcal Y}_2$ is isomorphic to the quotient of the commutative algebra  ${\mathcal R}[1,X,V]$ by the ideal generated by
$V^2-(d^2-1)V,VX-(d-d^{-1})V$ and $X^2=(d-2d^{-1})X+d^{-2}V$.}\\

\textbf{Theorem 1.2.} {\sl The algebra ${\mathcal Y}_3$ is generated by the 6 elements $A,B,C,F,G$ and $N$  in Figure 3.}\\
\begin{center}
\includegraphics[width=9cm,height=6cm]{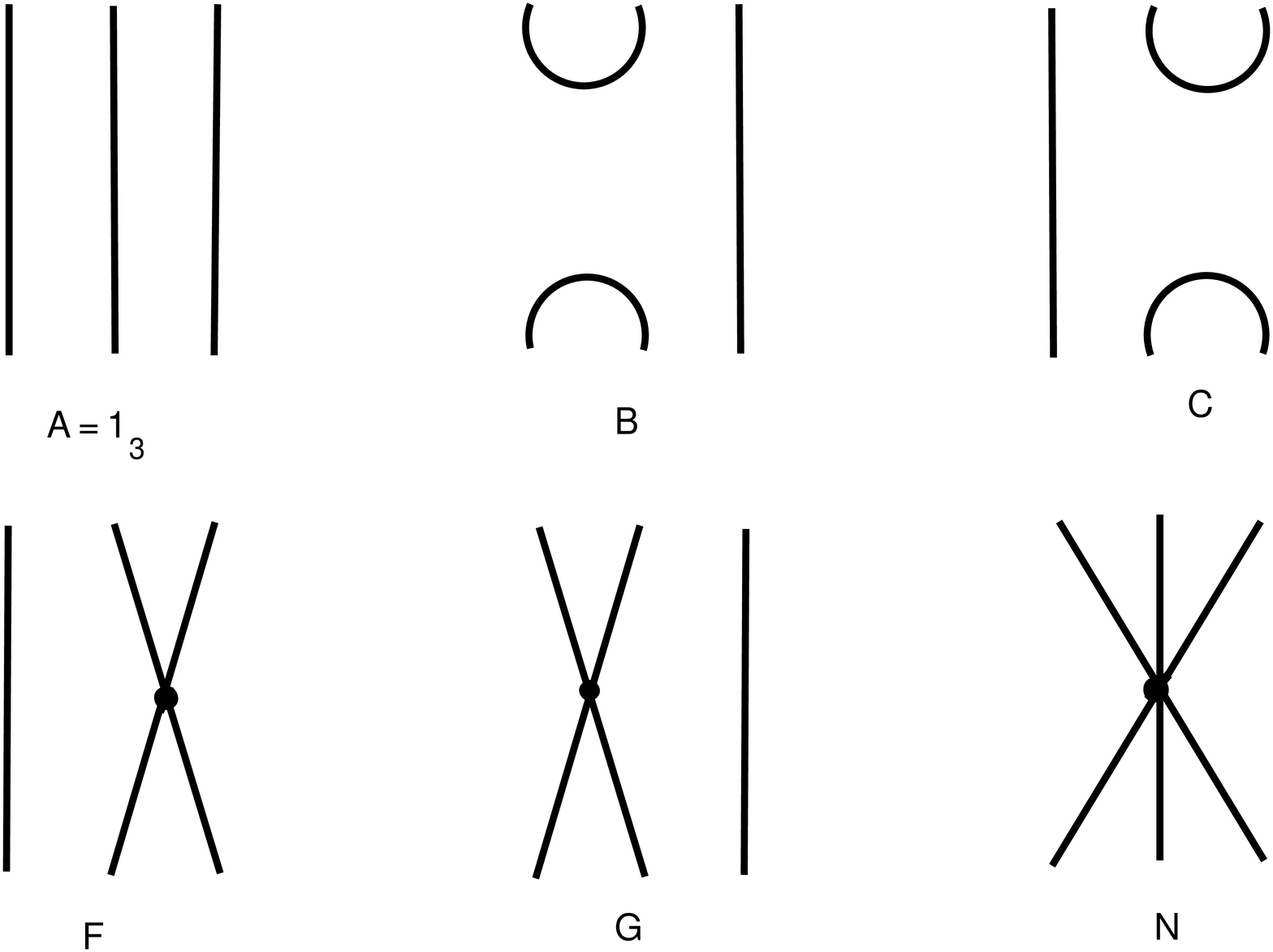} \\
Figure 3
\end{center}
This paper is outlined as follows. In Section 2, we briefly review some properties of the Temperley-Lieb algebra needed in the sequel. In Section 3, we define ribbon graphs and we introduce  the algebra ${\mathcal Y}_n$. The proofs of the main results  are given in  Section 4. Finally, Section 5 explores the connection between ${\mathcal Y}_n$   and the Temperley-Lieb algebra.
\section{The Temperley-Lieb algebra}
Temperley-Lieb algebras appeared first in the context of statistical physics. With the discovery of the Jones polynomial,  these algebras
  offered a new approach for the study of the   quantum invariants of links and three-manifolds. This section is  a brief introduction  to the theory of  Temperley-Lieb algebras from the knot theory viewpoint.

Let $n$ be a nonnegative  integer, an $n$-tangle $T$ is a one-dimensional
sub-manifold of $\R^2\times I$, such that the boundary  of $T$ is
made up of 2$n$ points $\{(i,0,0),(i,0,1);\;1\leq i\leq n\}$. As usual, tangles are considered up to isotopies  of  $\R^2\times I$ fixing the boundary pointwise. It is well known that the study of  tangles up to isotopy  is equivalent to the study of their planar diagrams in $\R \times I$ up to  Reidemeister moves keeping the boundary fixed pointwise.\\
Let
${\mathcal T}_n$ be the free $\mathcal R$-module generated by the
set of all $n$-tangles. We define $\tau _n$ to be the quotient of
${\mathcal T}_n$ by the smallest submodule containing all elements of the form:
 \begin{center}
$\hspace{2cm}  \bigcirc \cup
L - d  L $\\
$\hspace{2cm}  L -A L_{0} -A^{-1}  L_{\infty} ,$\\
\end{center}

 where $\bigcirc$ is the trivial circle, $d=-A^2-A^{-2}$ and  $L$, $L_0$ and $L_{\infty}$ are three tangle diagrams  which are identical everywhere  except in a small  disc where they look as pictured below.\\

\begin{center}
\includegraphics[width=10cm,height=2cm]{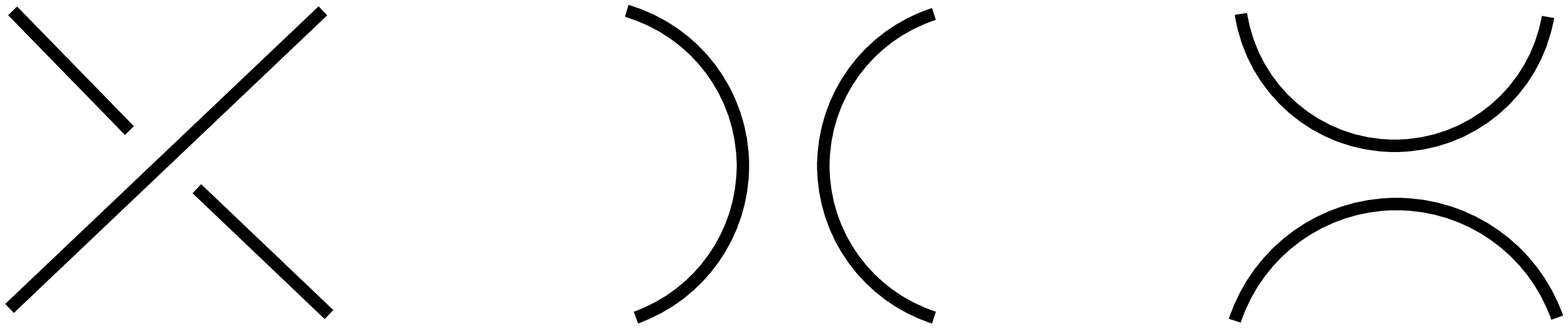} \\
{$L$}\hspace{3.5cm}{$L_0$}\hspace{3.5cm}$L_{\infty}$

\end{center}
If we equip  the module  $\tau _n$ with the standard  product of tangles, then $\tau _n$ turns out to be an algebra which is  isomorphic to the Temperley-Lieb algebra. A set of  generators
$(U_i)_{0\leq i \leq n-1}$ of $\tau _n$ is illustrated below.\\

\begin{center}

\includegraphics[width=10cm,height=4cm]{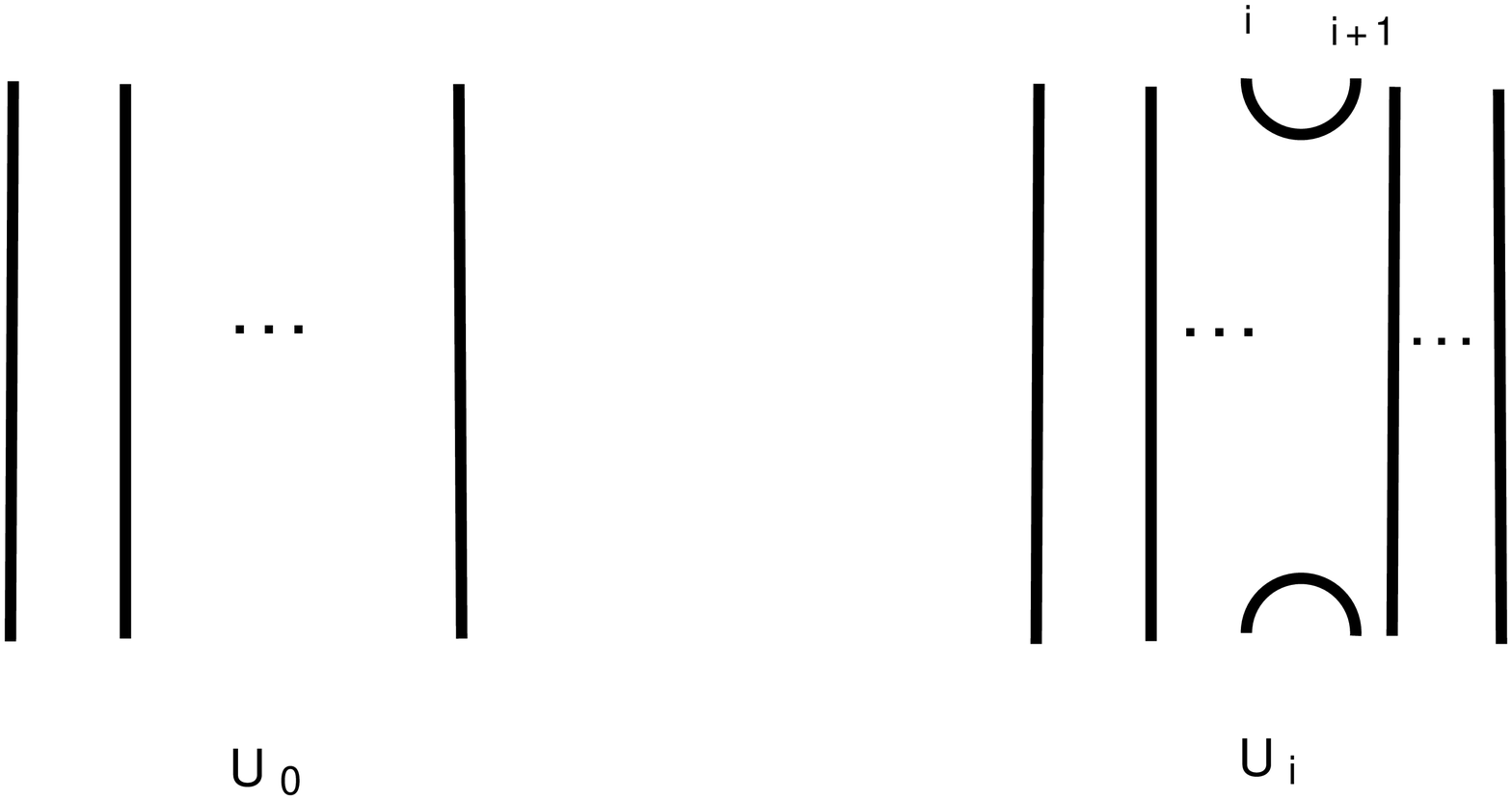}

  \end{center}

 Let $(f_i)_{0\leq i\leq n-1}$ denote the family of
Jones-Wenzl projectors in $\tau _n$. This family is  defined by
the
following recursive formulas:\\
$f_0=U_0,$\\
$f_{k+1}=f_k-\mu _{k+1}f_kU_{k-1}f_k$,\\
where $\mu _1=d^{-1}$ and $\mu _{k+1}=(d-\mu _k)^{-1}$.\\
In particular, we have  $f_1$=$1_n-d^{-1}U_1$. The elements
$f_k$ enjoy the following properties: $f_{k}^2=f_k$ and
$f_iU_j=U_jf_i=0$ for
$j\leq i$. See \cite{KL} for more details.
\section{Graph Algebra}
Ribbon graphs have been introduced by Reshetikhin and Turaev in \cite{RT}. They appeared  as a natural generalization of framed links.  We begin this section by a brief  review of the definition of ribbon graphs. More details about these objects can be found in \cite{Tu1}. Then, we will define the graph algebra $\mathcal{Y}_n$, which will appear as an extension of the theory of Temperley-Lieb algebra.\\
 Let $n \geq 1$ be an integer. A ribbon $n-$graph $G$  is a compact  oriented  surface embedded  into $\R^2\times I$ which can be decomposed  into  a finite  collection of annuli, coupons (small rectangles)  and ribbons (long bands), such that:\\
  (i) annuli do not meet each other and do not meet ribbons or coupons.\\
  (ii) ribbons  never meet each other, but they may meet coupons  at their bases.\\
  (iii) $G$ meets $\R^2\times \{0,1\}$ orthogonally exactly in bases of certain ribbons. The intersection is a collection of segments  $\{[i-1/10,i+1/10]\times\{0\}\times \{0,1\},i=1,\dots, n\}$.\\
   Ribbon $n-$graphs are considered up to isotopies  of  $\R^2\times I$ fixing the  boundary pointwise and preserving the decomposition into annuli, coupons and ribbons. According to \cite{Tu1}, ribbon $n-$graphs can be represented by
 planar diagrams, where coupons are represented by vertices, annuli are represented by circles  and long bands by either an ordinary graph edge, a half edge (arc connecting a vertex and a boundary point) or an arc connecting two boundary points.
\begin{center}
\includegraphics[width=4cm,height=3cm]{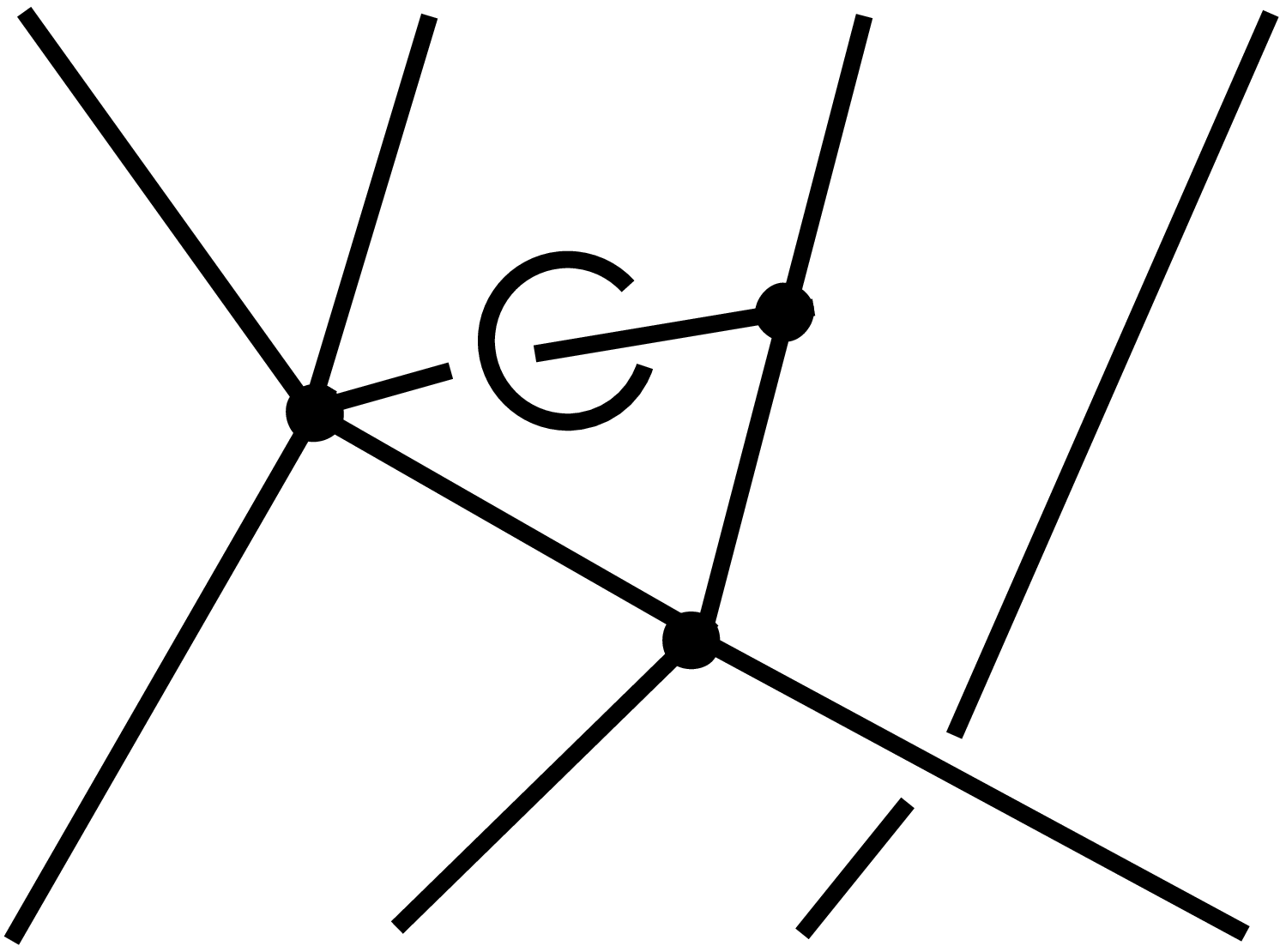}
\begin{center} {\sc Figure 4} \end{center}
 \end{center}

The study of ribbon graphs up to isotopy is equivalent to the study of their planar diagrams  modulo planar isotopies and   the  extended Reidemeister moves in Figure 5.\\

 \begin{center}
\includegraphics[width=4cm,height=1.5cm]{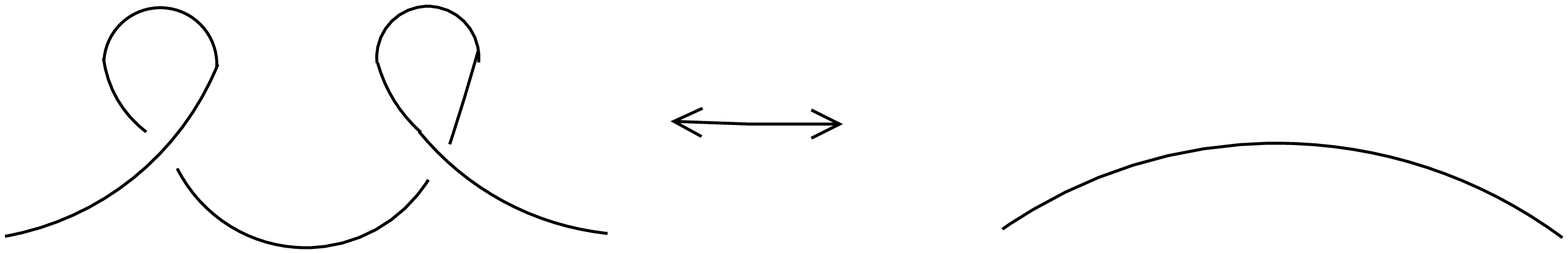}\\
\end{center}
\begin{center}\
\includegraphics[width=4cm,height=1.5cm]{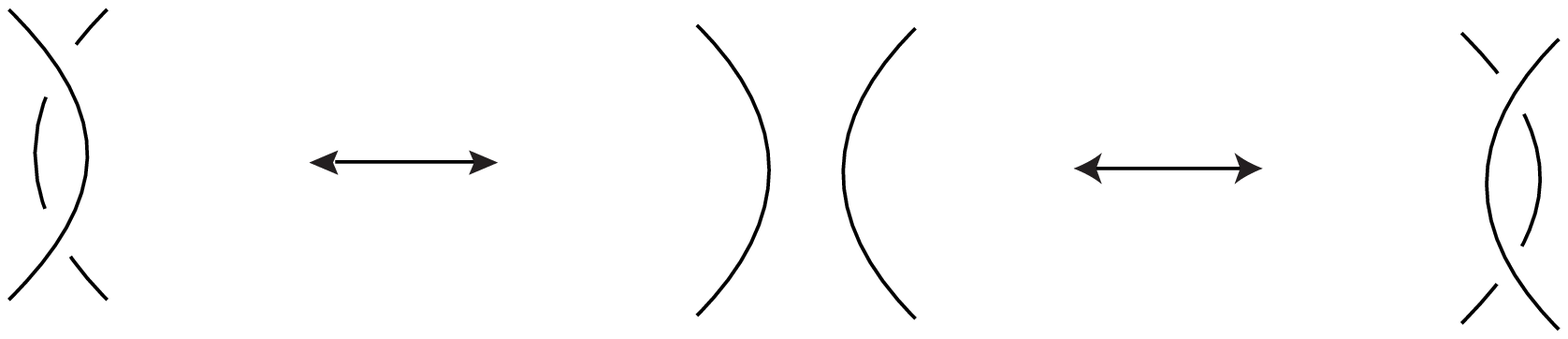}\\
\end{center}
\begin{center}
\includegraphics[width=4cm,height=1.5cm]{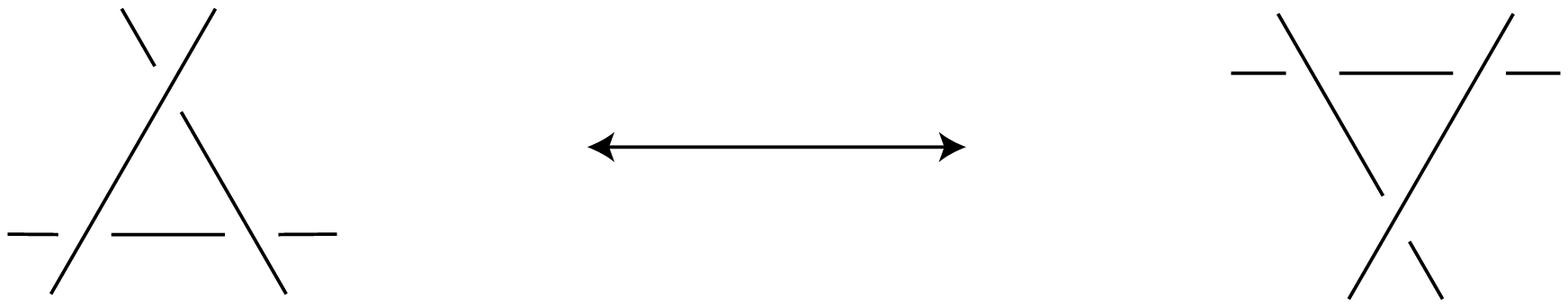}\\
\end{center}
\begin{center}
\includegraphics[width=4cm,height=1.5cm]{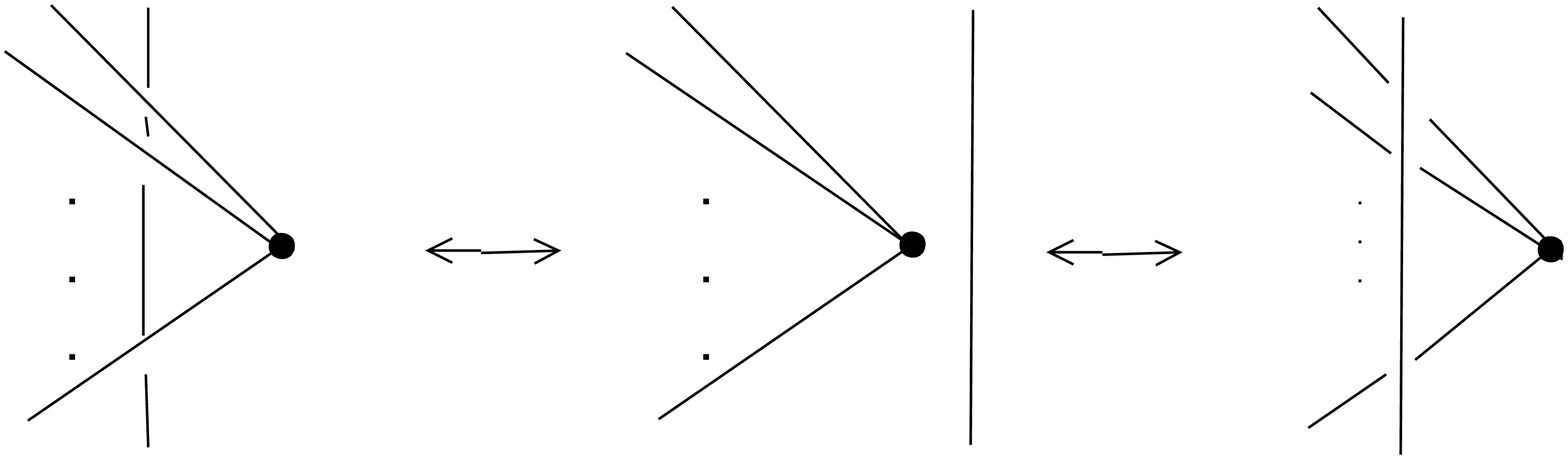}
\end{center}
\begin{center} {\sc Figure 5 } \end{center}

 Let ${\mathcal S}_n$ be the set of all ribbon $n-$graph diagrams and  $\mathcal
R$=$\Z[A^{\pm 1},d^{-1}]$, where $d=-A^2-A^{-2}$. Let
${Y}_n$ be the free $\mathcal R$-module generated by all elements
of  ${\mathcal S}_n$. We define $\mathcal{Y}_n$ to  be the quotient module  of ${Y}_n$ by the smallest submodule containing all expressions of the form:

\begin{center}
\includegraphics[width=12cm,height=9cm]{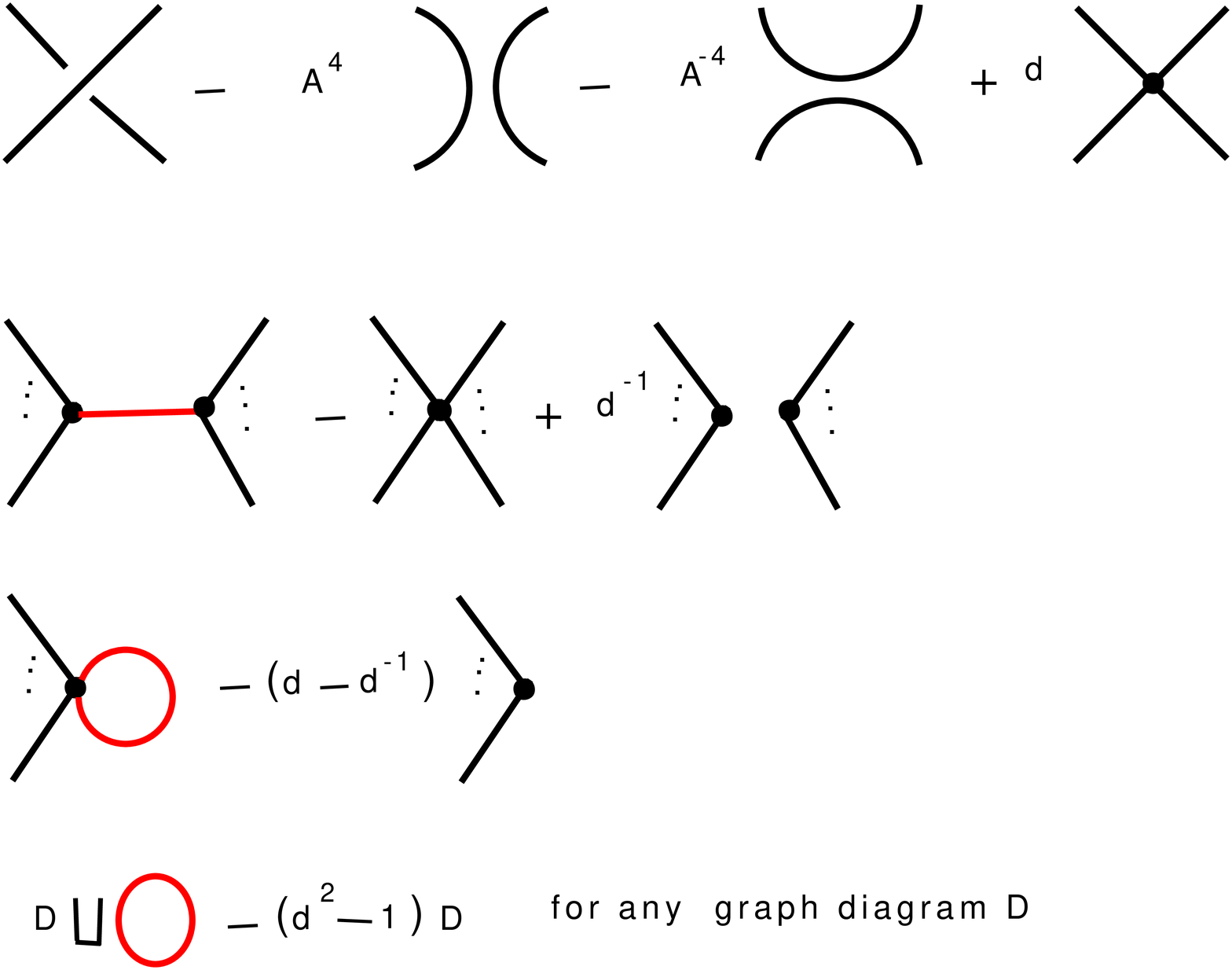} \\
\end{center}
 These relations are referred to as the Yamada skein relations. In each relation,  the pictures represent planar diagrams of ribbon $n-$graphs which are identical except in small disk where they look as pictured.
 In the same way as for tangles,  a multiplicative structure  can be defined on $\mathcal{Y}_n$. The identity relative to  this product is the ribbon graph  made up of $n-$parallel ribbons, this element is denoted hereafter by $1_n$. The product of two elements $G$ and $G'$ is the ribbon $n-$graph $GG'$ obtained by putting $G$ over $G'$ as pictured below:
 \begin{center}
\includegraphics[width=1.5cm,height=3cm]{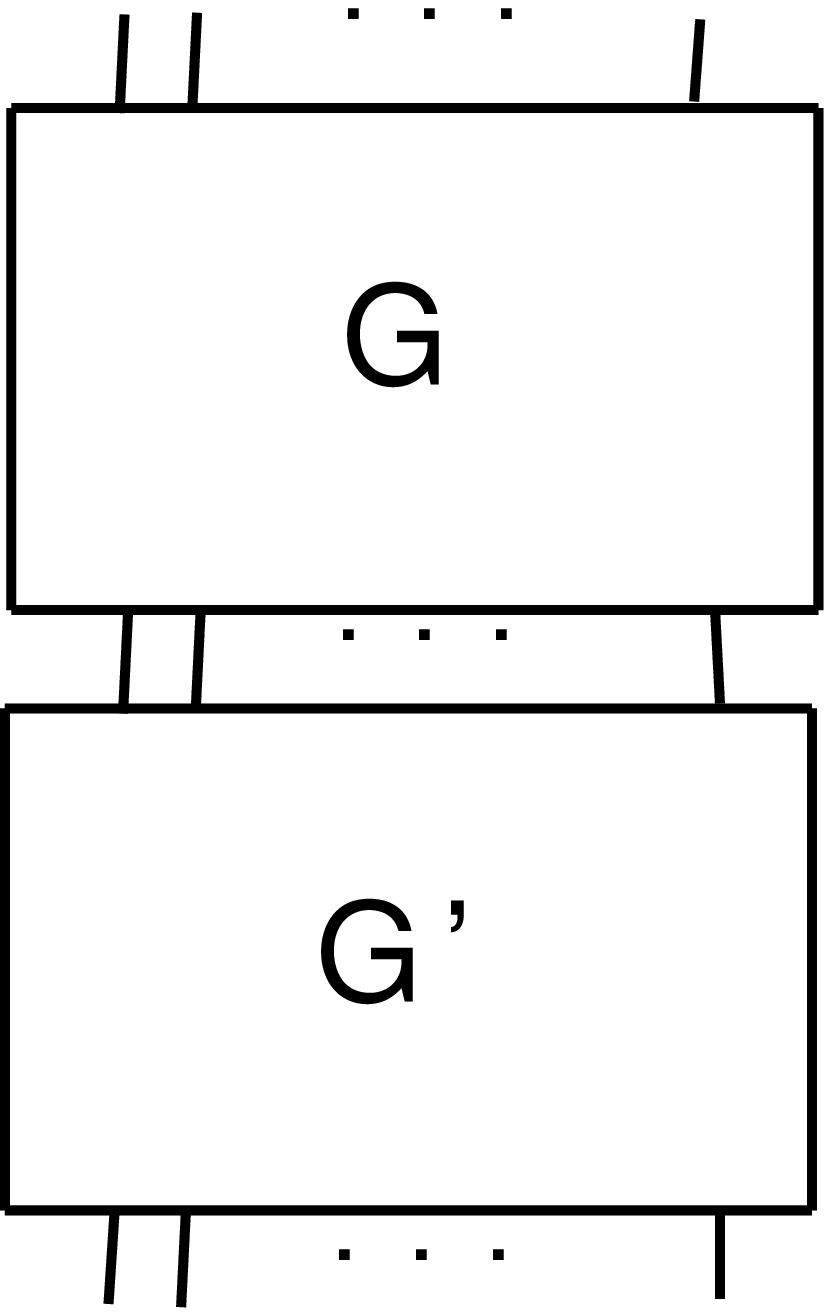}
\end{center}
\section{Proofs}
In this section we give the proofs of Theorem 1.1, Corollary 1.1 and Theorem 1.2. We begin by describing a family of generators of
 the ${\mathcal R}-$module ${\mathcal Y}_n$.  \\
\textbf{Lemma 4.1.}  \emph{The ${\mathcal R}-$module ${\mathcal Y}_n$ is generated by all ribbon $n-$graph diagrams with no crossings, no cycles  and no ordinary edges.}\\

\textbf{Proof.} Let $G$ be a ribbon $n$-graph diagram. One can apply  the first Yamada relation to smooth all the crossings of the diagram. Therefore, $G$ is  expressed  as a linear combination, with coefficients in $\mathcal R$ of ribbon $n-$graph diagrams each of which has no crossings. In the next step, we use the Yamada deletion-contraction relation  to delete all graph edges. Hence, our graph is written as a linear combination of diagrams which have  no ordinary edges. Now, we can remove all cycles  using Yamada relations (3) and (4). Finally, our graph $G$ is expressed as a linear combination of diagrams each of which has no crossings, no edges and  no cycles.\\
  \textbf{{Proof of Theorem 1.1.}}  According to  Lemma 4.1, the module ${\mathcal Y}_2$ is generated by the three elements
   $1_2$, $V$ and $X$ pictured in Figure 2.  So is the algebra  ${\mathcal Y}_2$.\\
  \textbf{Proof of Corollary 1.1.} The proof is straightforward by applying Yamada relations as illustrated below:\\
 $\begin{array}{lllll}
 &&& V^2=(d^2-1)V & \\
 &&&&\includegraphics[width=4.5cm,height=3cm]{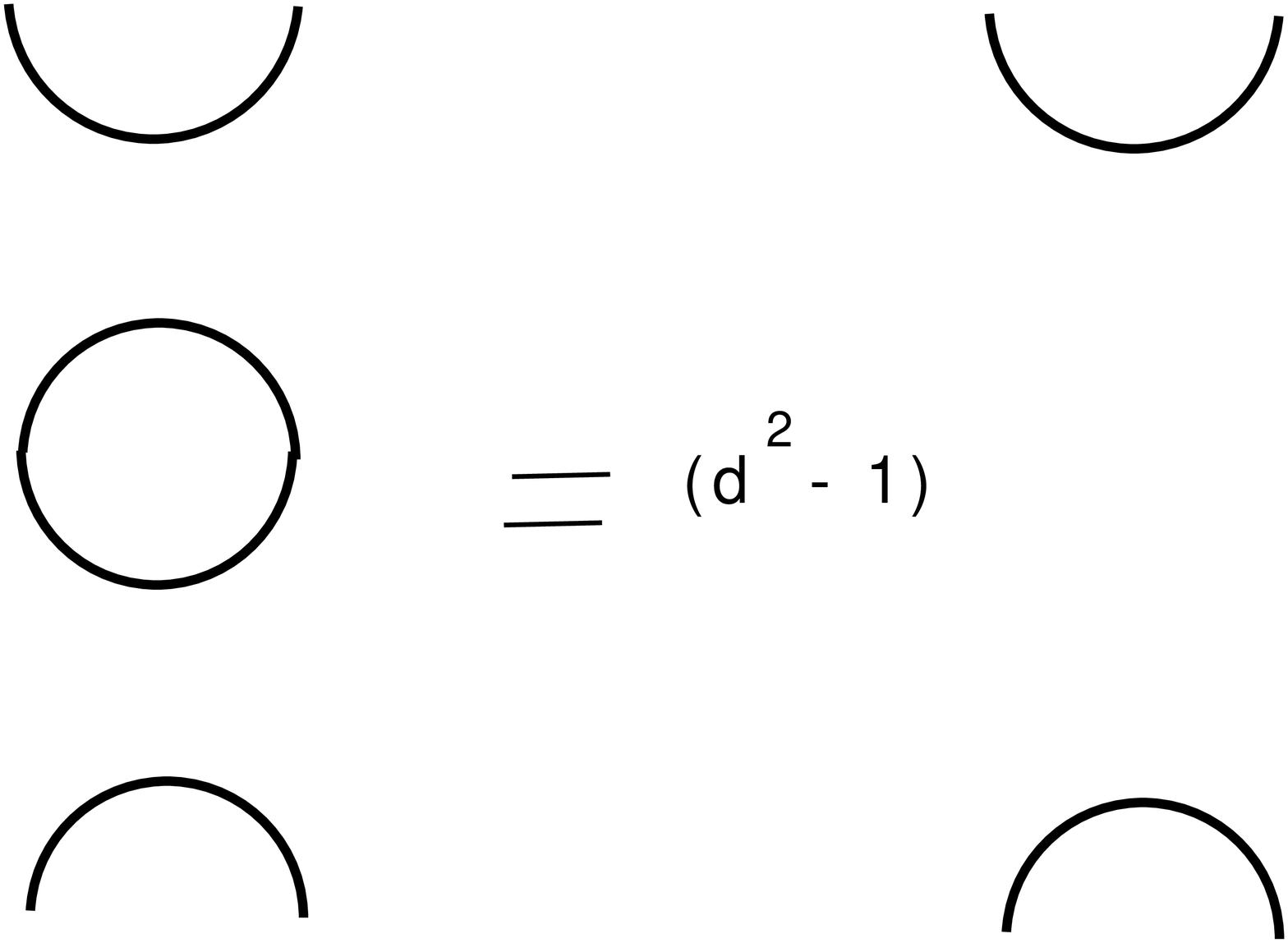}  \\
 && & VX=XV=(d-d^{-1})V &\\
  && & &\includegraphics[width=4.5cm,height=3cm]{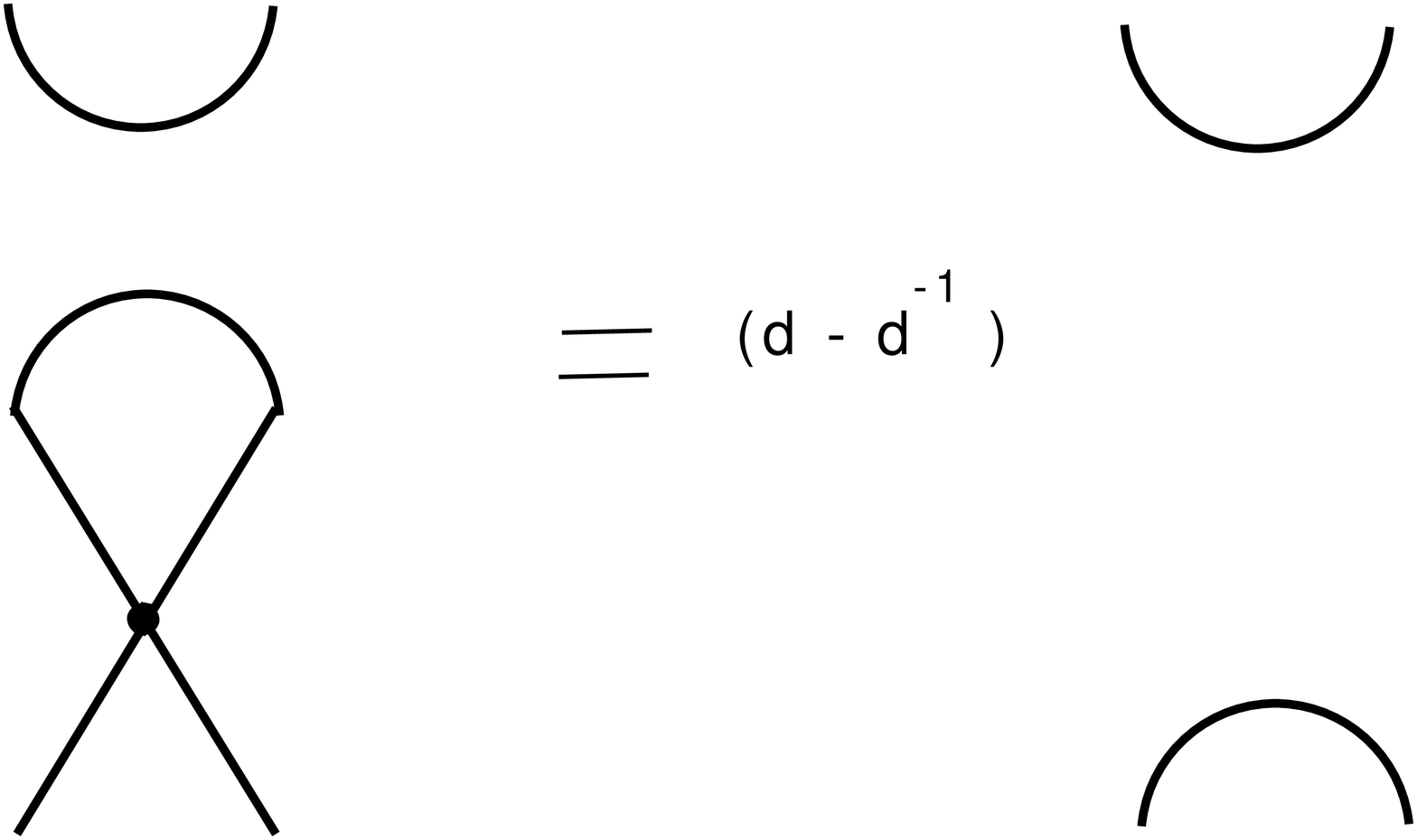} \\
  & &&X^2=(d-2d^{-1})X+d^{-2}V &\\
   && & &\includegraphics[width=8cm,height=5cm]{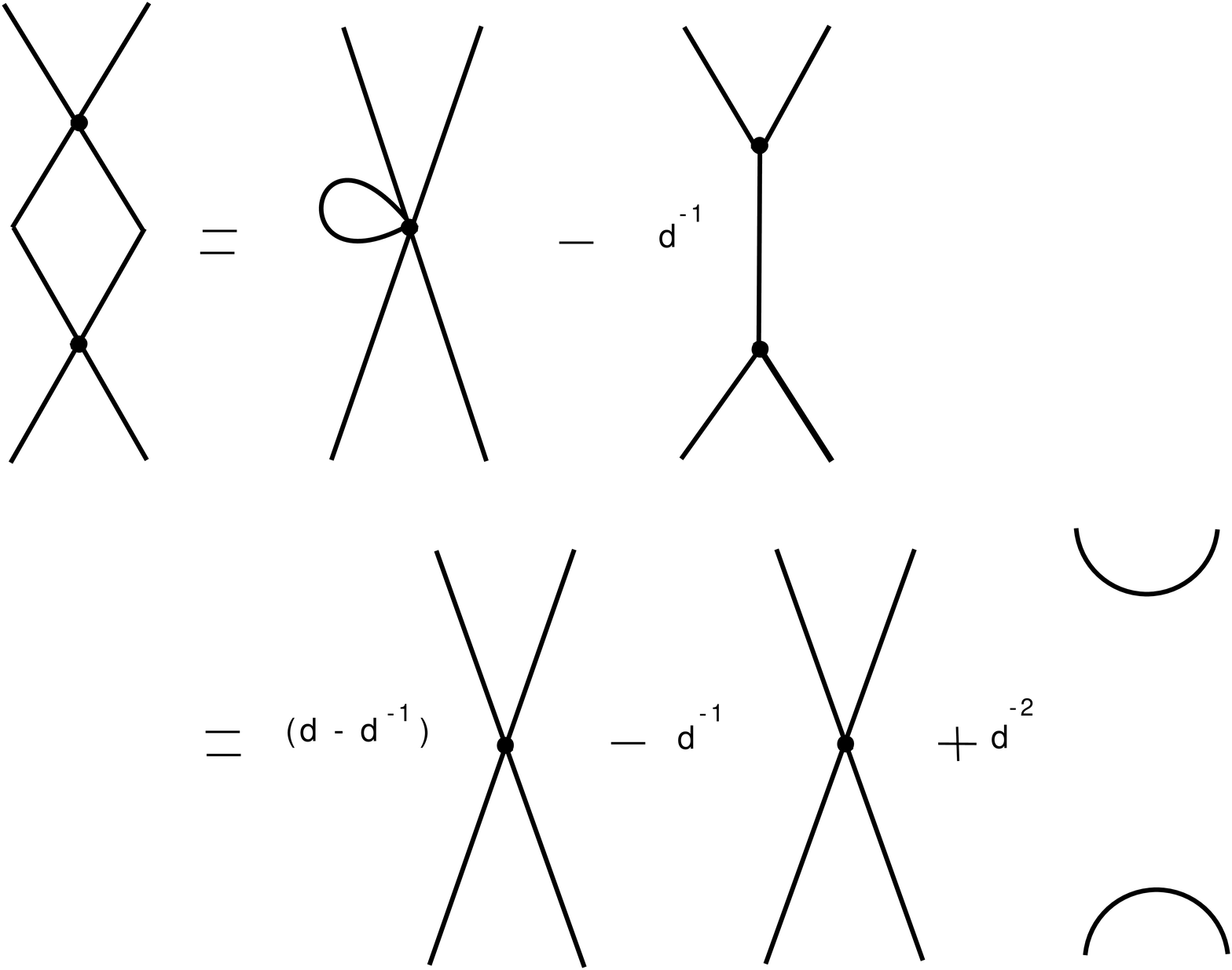} \\
  \end{array}$\\

\textbf{{Proof of Theorem 1.2.}} According to Lemma 4.1, the module   ${\mathcal Y}_3$ is generated by the following 15 elements:\\
  \begin{center}
\includegraphics[width=11cm,height=9cm]{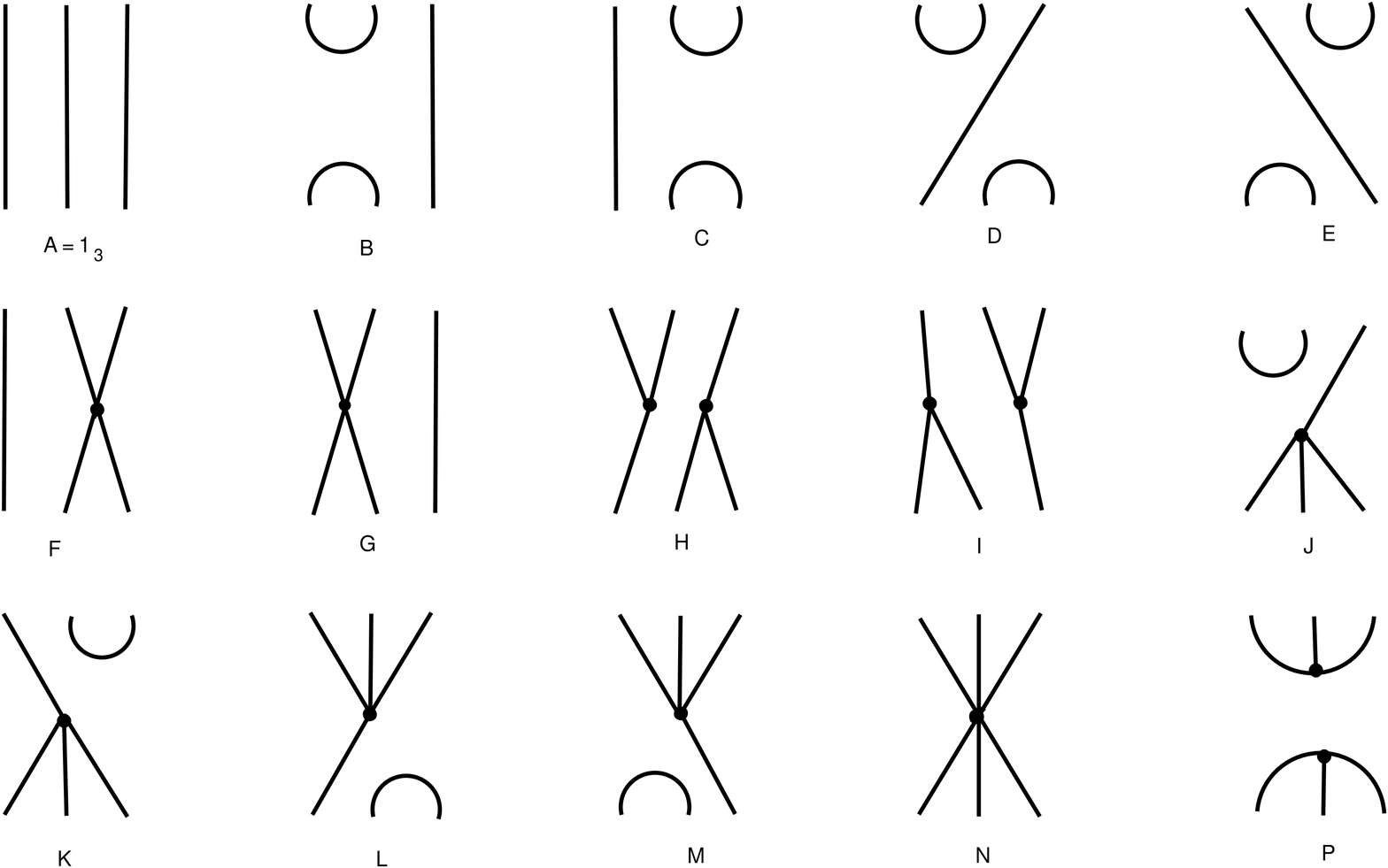} \\
\end{center}

  Now, we use the multiplication structure  to reduce these 15 generators of the module to  the six generators  of the  algebra ${\mathcal Y}_3$.
This reduction is briefly illustrated by  these 4 kind of operations\\
  1) It can be easily seen that: $D=BC$ and $E=CB$.\\
2) We have  $J=BF$ as depicted below
\begin{center}
\includegraphics[width=4cm,height=2cm]{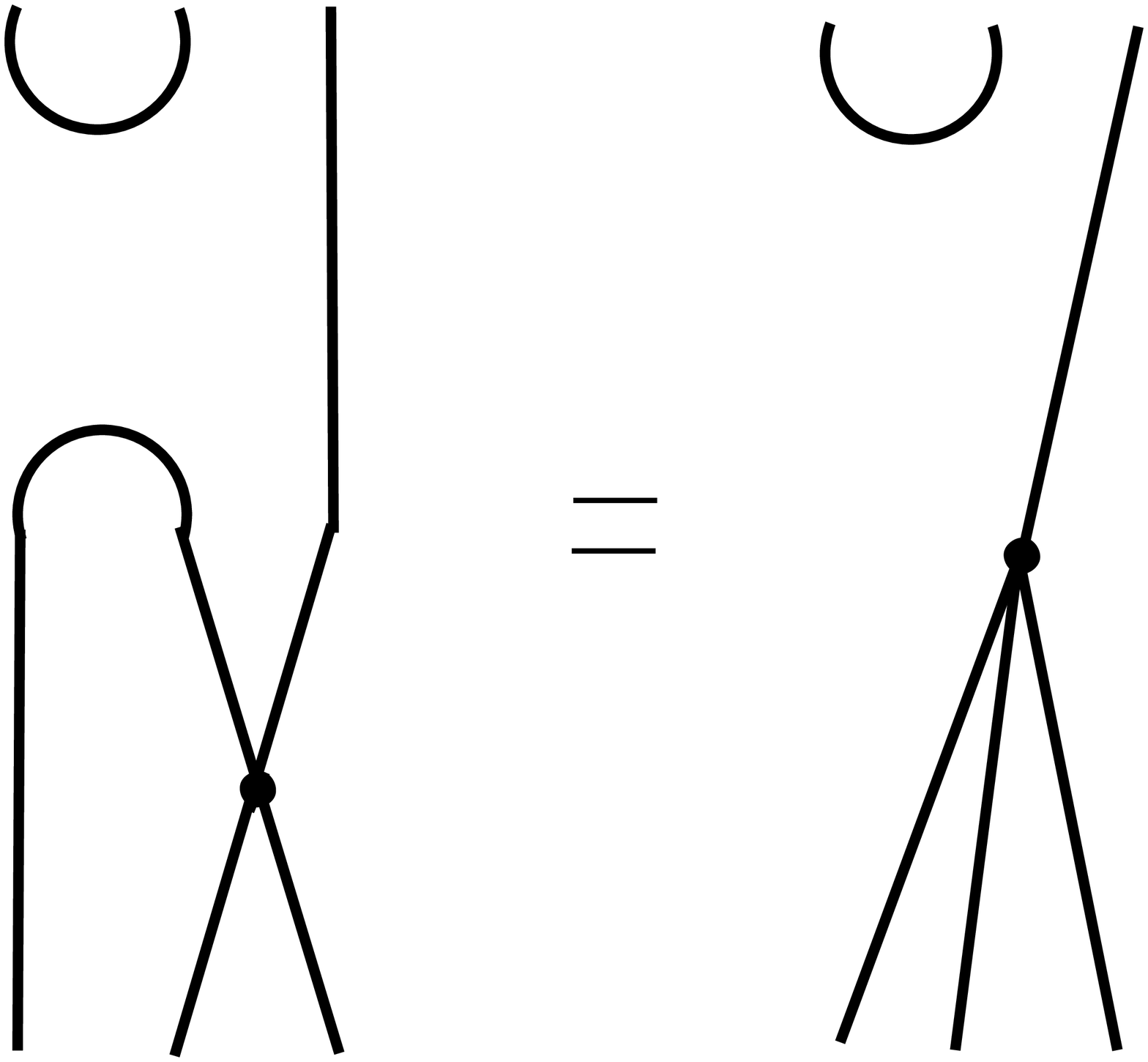} \\
\end{center}
Similarly we get  $M=FB$,  $K=CG$ and  $L=GC$.\\
 3)  According to the picture below we have: $FG=N-d^{-1}H$ which implies that $H=N-dFG$.\\
  \begin{center}
\includegraphics[width=8cm,height=2cm]{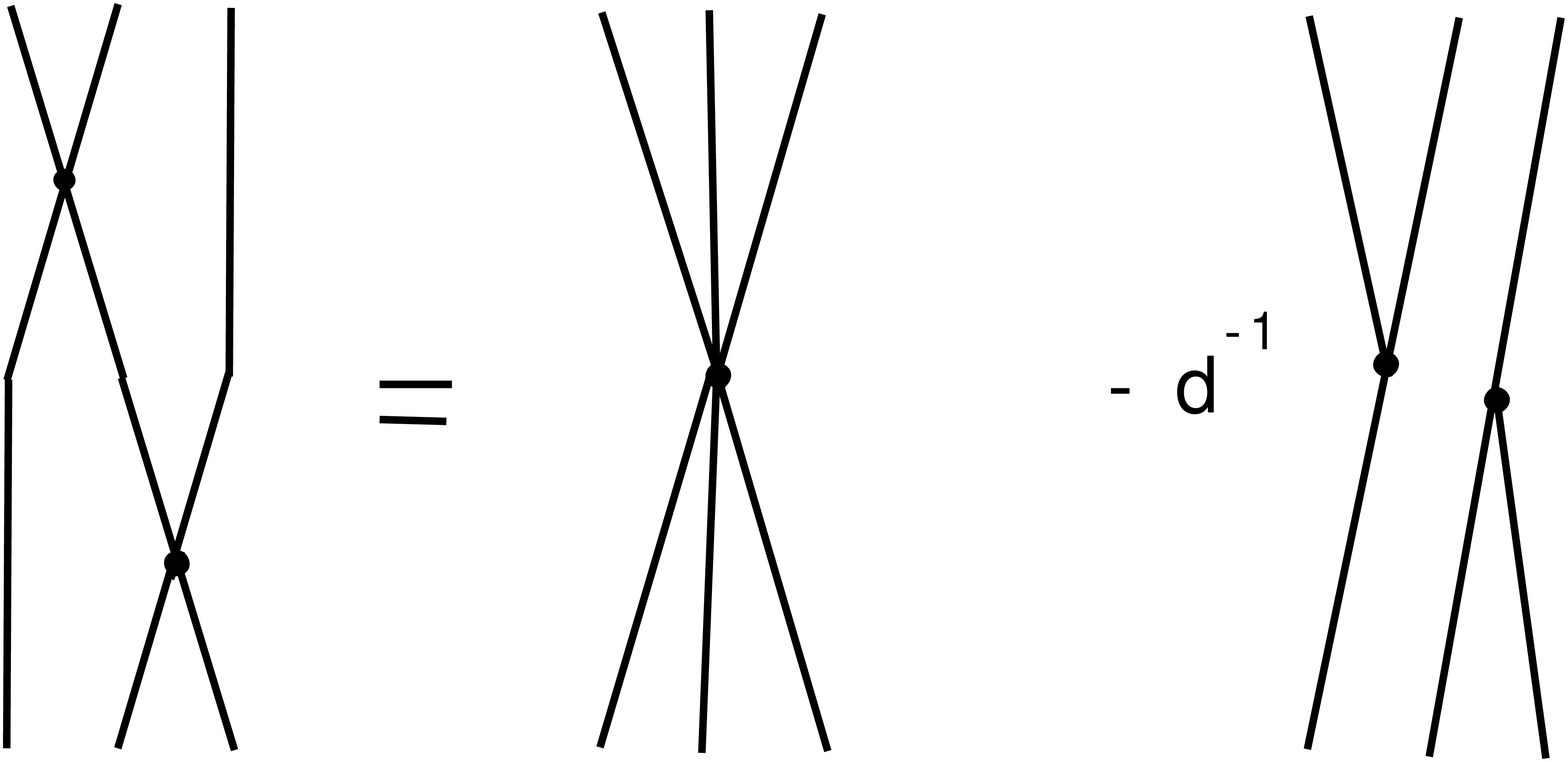} \\
\end{center}
Similarly, $GF=N-d^{-1}I$ which implies that $I=N-dFG$.\\
4) Finally,  using the deletion contraction formula as in the picture below,  we show that $MK=N-d^{-1}P$. This implies that:  $P=d(N-MK)=d(N-FBCG)$.
  \begin{center}
\includegraphics[width=8cm,height=2cm]{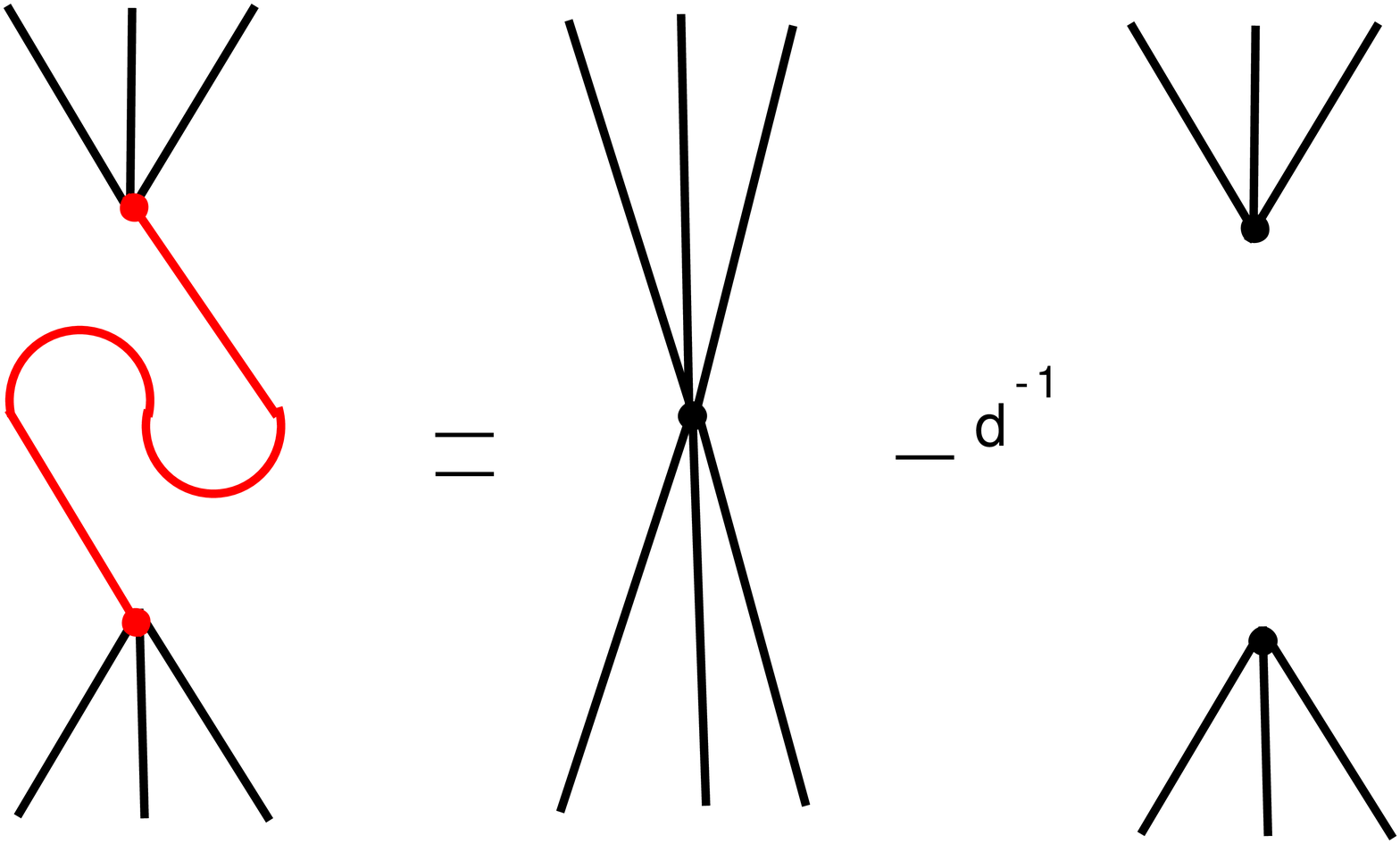} \\
\end{center}

\section{Relationship between ${\mathcal Y}_n$ and $\tau _n$} 
The purpose of this section  is to discuss  the relationship between the two  algebras ${\mathcal Y}_n$ and  $\tau _n$. In the case of skein modules of three-manifolds,  we   defined a  homomorphism  from the graph skein module to the Kauffman bracket skein module \cite{Ch1,Ch2}. An analogous s of this homomorphism is defined here. Let ${\varphi}_n: Y_n \longmapsto \tau_{2n} $ be the linear map that associates to each ribbon $n$-graph diagram $G$
the  linear combination of  diagrams  obtained from $G$ by replacing each edge, half edge and arc  of $G$ by two
planar strands with a projector $\begin{picture}(0,0)
\put(0,0){$f_{1}=$} \put(30,-10){\line(0,1){25}}
\put(37,-10){\line(0,1){25}} \put(40,0){$-d^{-1}$}
\put(75,14){\oval(7,17)[b]} \put(75,-10){\oval(7,17)[t]}
\end{picture}$\\ in the cable, and by replacing each vertex of $G$ by a diagram as follows (the figure illustrates  the case of a four-valent vertex)\\
\begin{center}
\includegraphics[width=8cm,height=3cm]{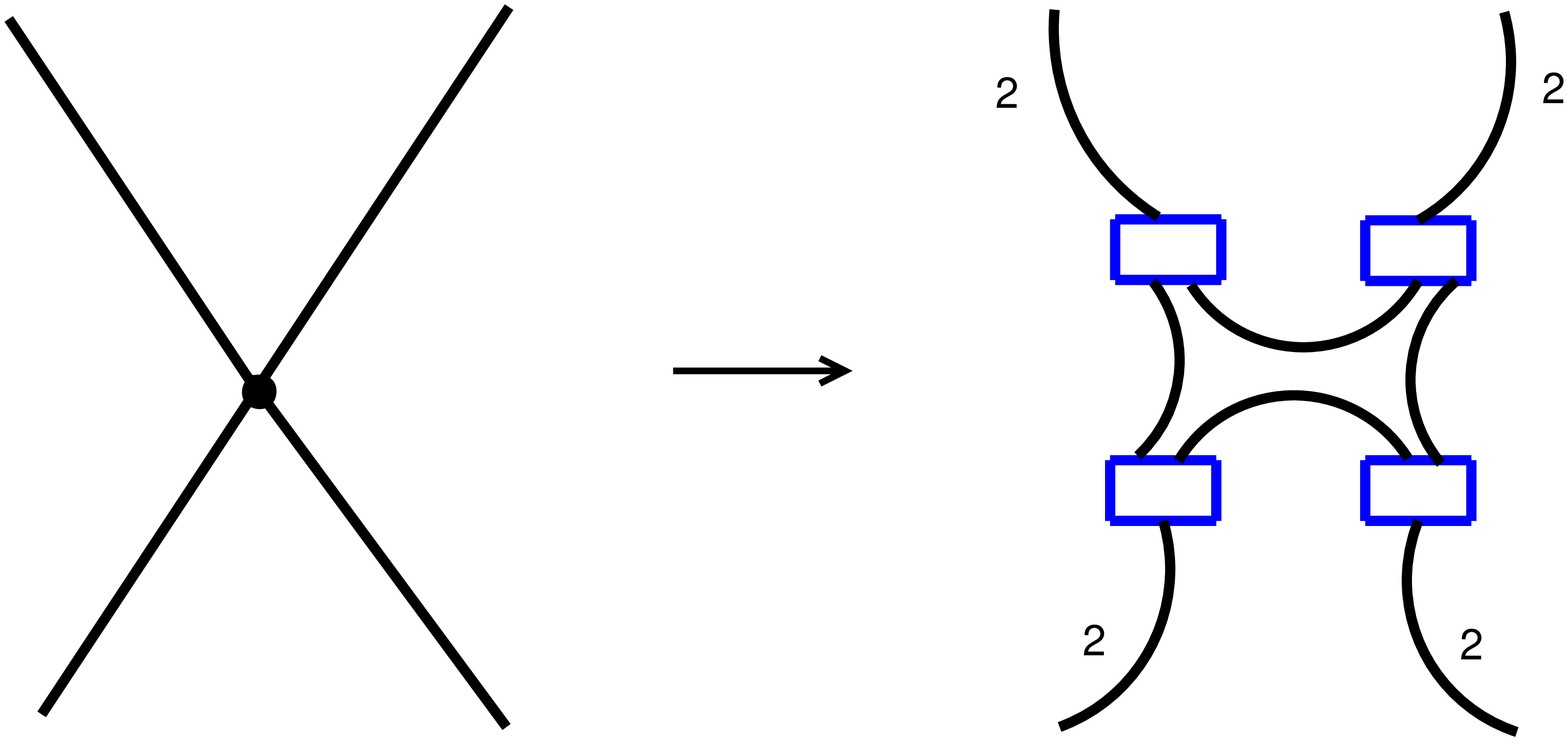}
\end{center}
\begin{center} {\sc  Figure 6 } \end{center}

In this picture, writing an integer $2$ beneath an edge  $e$ means that this
edge has to be replaced by $2$ parallel ones. Notice that  ${\varphi}_n$ is defined on the generators of the free $\mathcal R$-module, then extended by linearity to all elements of $Y_n$. Using the same graphic calculations as in \cite{Ch1} (Lemma 3.4), one can check easily that  $\varphi_n$ defines a map $\Phi_n$ from  ${\mathcal Y}_n$ to $\tau_{2n}$.  Obviously,  $\Phi_n$ is a homomorphism of algebras. The following picture illustrates how to compute  $\Phi_2(1_2)$. \\
\begin{center}
\includegraphics[width=10cm,height=2cm]{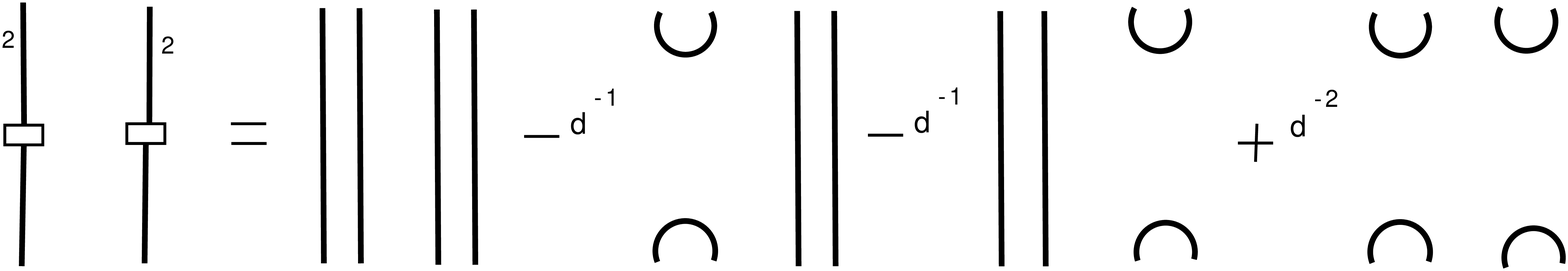}
\end{center}

\textbf{Theorem 5.1.} The homomorphism  \emph{ $\Phi_n: {\mathcal Y}_n \longmapsto \tau_{2n} $ is injective.} \\
\textbf{Proof.}  Remind first that the $\mathcal R$-module $\tau_n$ has a standard base consisting of all diagrams of $n-$arcs with no crossings joining the $2n$ boundary points pairwise. The dimension of $\tau_n$ is   the Catalan number $\displaystyle\frac{C_n^{2n}}{n+1}$. Now, let $t$ be an element of the standard  base of $\tau_{2n}$. If we take the union of $t$ with the $n$ segments $[i,i+1]\times \{0\} \times \{0,1\}$ for $i$ odd, then we get a 1-dimensional manifold which bounds a surface $\Sigma(t)$ in $\R \times [0,1]$. The surface  $\Sigma(t)$ retracts by deformation (in  $\R \times [0,1]$ ) on an  graph diagram $g(t)$.  The picture below  illustrates this construction in the case of an element  $ t\in \tau_6$.\\

\begin{center}
\includegraphics[width=11cm,height=3cm]{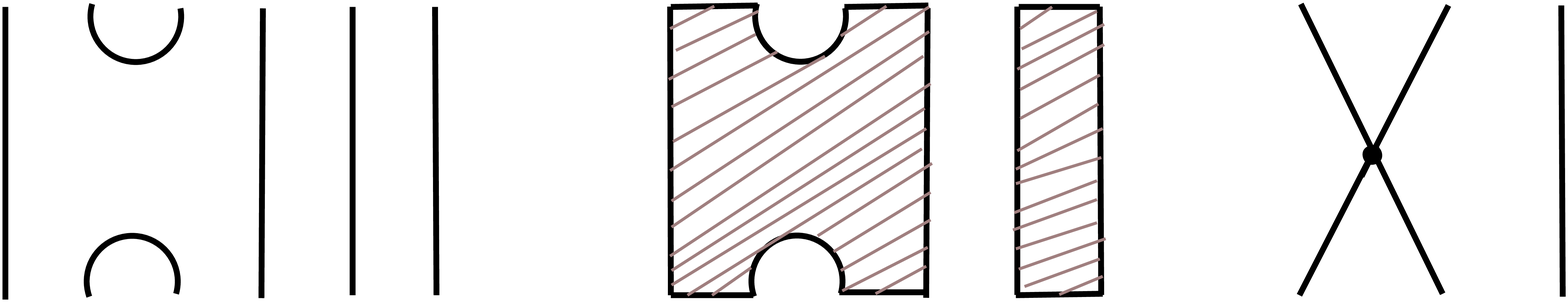}\\
{$t$}\hspace{3.5cm}{$\Sigma(t)$}\hspace{3.5cm}$g(t)$
\end{center}
\begin{center} {\sc  Figure  7} \end{center}

Now, we shall   prove that the kernel of $\Phi_n$ is trivial. Let $g_1, \dots, g_s$ be distinct generators  of  ${\mathcal Y}_n$  as described in Lemma 4.1. Let  $r_1, \dots, r_s$ be elements of $\mathcal R$ such that $\Phi_n(r_1g_1+\dots+r_sg_s)=0$.\\
We know that  $\Phi_n(g_i)$ is expressed as a linear combination of  the standard generators  of $\tau_{2n}$.  Among the elements which appear in this combination, let $g_i^1$ be the generator  whose surface $\Sigma (g_i^1)$ has the minimum connected components. Obviously,  this surface retracts by deformation on $g_i$. Moreover, it   is easy to see that  $\Phi_n(r_1g_1+\dots+r_sg_s)=0$ implies that $r_1g_1^1+\dots+r_sg_s^1=0$ which leads to $r_i=0$ for all $1\leq i \leq s$. Hence, $\Phi_n$ is injective.\\


\begin{thebibliography}{99}

 \bibitem{Ch1} N. {\sc Chbili}.
{\em Skein algebras of the solid torus and symmetric spatial graphs}. Fund. Math. 190, (2006), pp. 1-10.
\bibitem{Ch2} N. {\sc Chbili}.
{\em Graph skein modules and symmetries
of spatial graphs}. Preprint.
\bibitem{CF} N. {\sc Chbili} and T. {Fleming}. {\em The graph skein algebras of the torus}. In preparation.
\bibitem{Ka} L. {\sc Kauffman}. {\em An invariant of regular isotopy.} Trans. Amer. Math.
Soc. 318 (1990), no. 2, pp. 417-471.
\bibitem{KL} L. H. {\sc Kauffman} and S. L {\sc Lins}. {\em Temperley-Lieb recoupling theory and invariants of 3-manifolds.}
 Ann. Math. Studies. 134, Princeton Univercity Press (1994).

\bibitem{Li} W. B. R.  {\sc Lickorish}. {\em The skein method for 3-manifold invariants.}
J. Knot Th. Ram. 2 (1993), 171-194.

 \bibitem{Pr1} J. H.  {\sc Przytycki}. {\em  Skein modules of 3-manifolds.} Bull. Pol. Acad. Sci.: Math.,
39, 1-2 (1991) pp. 91-100.

 \bibitem{Pr3} J. H.  {\sc Przytycki}.  {\em Fundamentals of Kauffman bracket skein modules.} Kobe Math. J., 16(1), (1999), 45-66.

\bibitem{RT}  N. {\sc Reshetikhin}  and V. G.  {\sc Turaev}. {\em Ribbon Graphs and Their Invariants Derived
from Quantum Groups.} Commun. Math. Phys. 12 (1990),  pp. 1-26

 \bibitem{Tu1} V. G.  {\sc Turaev}. {\em Quantum invariants of knots and 3-manifolds.} De Gruyter, 2nd edition 2010.
 \bibitem{Ya1} S. {\sc Yamada}.  {\em An invariant of  spatial
 graphs.}
 J. Graph theory, 13 (1989) pp. 537-551.
 \bibitem{Ya2} S. {\sc Yamada}  {\em A topological invariant of  spatial
 regular graphs.}
 Proceeding of Knots 90, De Gruyter 1992,  pp. 447-454.
 \bibitem{Yo} Y. {\sc Yokota}  {\em Topological invariants of
 graphs in 3-space.}
 Topology, Vol. 35, (1996), pp. 77-87.


\end{thebibliography}
\end{document}